\documentclass[journal,10pt,twocolumn,twoside]{IEEEtran}

\usepackage{amsmath,amssymb,amsthm,url,amsfonts,mathrsfs, bbm}
\usepackage{graphicx, epsfig}
\usepackage{verbatim}
\usepackage{subfigure}
\usepackage{ifthen}
\usepackage{multirow}
\usepackage{multicol}
\usepackage{float}
\usepackage{subeqnarray}
\usepackage{enumerate}
\usepackage{cite}
\usepackage[usenames,dvipsnames]{color}
\usepackage[font=normalsize]{caption}

\usepackage[bookmarks=true,pageanchor,colorlinks,linkcolor=blue,anchorcolor=blue,citecolor=blue,urlcolor=blue,hyperfootnotes=false]{hyperref}
\usepackage{color}

\usepackage{ifthen}
\newboolean{showcomments}
\setboolean{showcomments}{false}

\newcommand{\eyb}[1]{\ifthenelse{\boolean{showcomments}}
{ \textcolor{red}{(Eilyan says:  #1)}}{}}
\newcommand{\bose}[1]{  \ifthenelse{\boolean{showcomments}}
{ \textcolor{red}{(Bose says:  #1)}} {}  }
\newcommand{\rev}[1]{#1}

\newcommand{\beq}{\begin{equation}}
\newcommand{\eeq}{\end{equation}}
\newcommand{\beqa}{\begin{eqnarray}}
\newcommand{\eeqa}{\end{eqnarray}}
\newcommand{\beqan}{\begin{eqnarray*}}
\newcommand{\eeqan}{\end{eqnarray*}}

\newcommand{\E}{\mathbb{E} }

\newcommand{\prob}{\mathbb{P}}

\newcommand{\Rset}{\mathbb{R}}

\newcommand{\Acal}{{\cal A}}
\newcommand{\Bcal}{{\cal B}}
\newcommand{\Ccal}{{\cal C}}

\newcommand{\Fcal}{{\cal F}}

\newcommand{\Kcal}{{\cal K}}
\newcommand{\Lcal}{{\cal L}}

\newcommand{\Ncal}{{\cal N}}

\newcommand{\Pcal}{{\cal P}}

\newcommand{\Scal}{{\cal S}}

\newcommand{\Zcal}{{\cal Z}}

\newcommand{\RR}{{\bf R}}

\newcommand{\bone}{\mathbf{1}}
\renewcommand{\v}[1]{{\mathbf{#1}}}
\newcommand{\ve}{\varepsilon}
\newcommand{\ol}[1]{\ensuremath{\overline{{#1}}}}

\newcounter{l1}
\newcounter{l2}
\newcounter{l3}
\newcommand{\bdotlist}{\begin{list}{$\bullet$}{}}
\newcommand{\bboxlist}{\begin{list}{$\Box$}{}}
\newcommand{\bbboxlist}{\begin{list}{\raisebox{.005in}{{\tiny
$\blacksquare$ \ \ }}}{}}
\newcommand{\bdashlist}{\begin{list}{$-$}{} }
\newcommand{\blist}{\begin{list}{}{} }
\newcommand{\barablist}{\begin{list}{\arabic{l1}}{\usecounter{l1}}}
\newcommand{\balphlist}{\begin{list}{(\alph{l2})}{\usecounter{l2}}}
\newcommand{\bAlphlist}{\begin{list}{\Alph{l2}.}{\usecounter{l2}}}
\newcommand{\bdiamlist}{\begin{list}{$\diamond$}{}}
\newcommand{\bromalist}{\begin{list}{(\roman{l3})}{\usecounter{l3}}}

\newtheorem{theorem}{Theorem}
\newtheorem{lemma}{Lemma}
\newtheorem{proposition}{Proposition}

\newtheorem{remark}{Remark}
\newtheorem{assumption}{Assumption}

\theoremstyle{definition}
\newtheorem{definition}{Definition}

\renewcommand{\Rset}{\mathbf{R}}

\renewcommand{\bone}{\mathbf{e}}
\renewcommand{\v}[1]{#1}

\newcommand{\vpi}{{\v{\pi}}}
\newcommand{\interior}[1]{\text{int}({#1})}
\newcommand{\bound}[1]{\partial{#1}}

\newtheorem*{lemma*}{Lemma}

\title{
The Marginal Value of Networked Energy Storage
}

\author{Subhonmesh Bose$^1$,~\IEEEmembership{Member,~IEEE,} and Eilyan Bitar$^2$,~\IEEEmembership{Member,~IEEE,}
\thanks{This work was supported in part by NSF grant ECCS-1351621, NSF grant CNS 1239178, PSERC sub-award S-52, the US DoE under the CERTS initiative, and Cornell University's David R. Atkinson Center for a Sustainable Future (ACSF).}%
\thanks{$^1$ S. Bose is with the Department of Electrical and Computer Engineering, University of Illinois at Urbana-Champaign, Urbana, IL 61801, USA. ({\tt\small boses@illinois.edu})}%
\thanks{$^2$ E. Bitar  is with the School of Electrical and Computer Engineering, Cornell University, Ithaca, NY 14853, USA. ({\tt\small eyb5@cornell.edu})}%

}
\date{}

\begin{document}
\maketitle

\begin{abstract}
We consider the problem of characterizing the locational marginal value of energy storage capacity in electric power networks with stochastic  renewable supply and demand. The perspective taken is that of a system operator, whose objective is to minimize the expected cost of firm supply required to balance a stochastic net-demand process over a finite horizon, subject to transmission and energy storage constraints. The value of energy storage capacity is defined in terms of the optimal value of the corresponding constrained stochastic control problem.
It is shown to be concave and non-decreasing in the vector of location-dependent storage capacities --  implying that the greatest marginal value of storage is derived from initial investments in storage capacities. 
We also provide -- as part of our main result -- a characterization of said marginal value, which reveals its explicit dependency on a specific measure of nodal price variation. 
More generally, we derive an upper bound on the locational marginal value of energy storage capacity in terms of the total variation of the corresponding nodal price process, and provide conditions under which this bound is tight.

\end{abstract}


\section{Introduction}

The variability in supply inherent to renewable energy resources like wind and solar poses a fundamental challenge to their integration into power system operations at scale. The primary difficulty derives from the need to instantaneously balance an inelastic demand for power with an intermittent supply of power across a transmission-constrained network. Energy storage  offers a form of flexibility that enables the absorption of power imbalances through appropriate reshaping of demand and supply profiles over time. The ability to do so should in turn enable a substantial reduction in the cost of operating a power system with significant penetration of renewable energy resources. Naturally, this reduction in cost  depends critically on the collective placement, sizing, and control of energy storage assets. In order to gauge the potential impact of this emerging technology, one needs to accurately quantify the value of storage -- measured in terms of  reduction in system operating cost. This paper offers a mathematical framework to tractably quantify this value. Specifically, we provide an expression for the location-specific value derived from \emph{initial investment} in storage capacity over a constrained power network with stochastic  supply and demand. 
This value is shown to depend  on both properties of the power network and the stochastic processes driving the system. From hereon, we refer to \emph{net demand} as the  nominal demand minus the variable supply.

Many have advocated the adoption of storage technologies to enable the transition to a future power system with deep integration of renewables \cite{iecWhitePaper, eyer2010energy, rittershausen2011moving, denholmnrel}. Storage assets come in different modalities. Electrochemical batteries, flywheels, and pumped hydro are the most common examples. In addition,  aggregated flexible loads can provide storage services \cite{callaway2009tapping, taylor2013competitive, hao2013generalized}. Of interest are questions pertaining to the optimal placement, sizing, and control of such storage assets over a power network with stochastic demand and supply. As we explain through our formulation, such problems in their full generality are not conducive to tractable mathematical analyses. As a result, different papers resort to different sets of simplifying assumptions. We summarize a subset of the related literature here. Su and El Gamal \cite{su2011modeling} and Parandehgheibi et al. \cite{dahleh12} study the optimal storage control problem in a stochastic control framework. Their analysis utilizes the so-called copperplate model, where the power system is treated as a single bus network.
On the other extreme, Thrampoulidis et al. \cite{thrampoTAC}, Bose et al. \cite{Bose2012storage}, and Castillo and Gayme \cite{anyaCDC} study the joint problem of optimal placement and control of storage resources over a power network, but with a deterministic model of net demand. 
None of the aforementioned papers \cite{su2011modeling, dahleh12, thrampoTAC, Bose2012storage, anyaCDC} accomodate  both transmission constraints and uncertainty in net demand in their analyses.
An exception is the work by Kanoria et. al \cite{kanoria2011distributed}. 
Recognizing the difficulty inherent to the explicit treatment of linear transmission constraints, Kanoria et. al \cite{kanoria2011distributed}  augments the objective function to include quadratic penalties on violation of such constraints. Given the augmented problem, they are able to characterize  the optimal control policy for power networks with a regular topology and a stationary net demand process.

The above papers all model the storage assets as being controlled by a system operator (SO) with an aim to optimize a system-wide objective. A long list of papers have, in parallel, considered the problem of an individual storage owner-operator, who aims to maximize her expected revenue from energy sales in an electricity market. Examples include \cite{baldick, ramRaja, harshaoptimal, sioshansi, vandeven, KimPowell2011, bitar2011role,  bitar_trans, qin2013storage, thatte2013risk, raovalue}.

\subsubsection*{Our contribution} This paper considers the SO's problem of determining a system dispatch  policy to minimize the expected cost of balancing a stochastic net demand process across a transmission-constrained power network -- the so called multi-period economic dispatch problem. Formulating the SO's problem as a finite horizon, stochastic optimal control problem with perfect state feedback, we first establish in Theorem \ref{thm:convex}, the  
convexity and monotonicity of the optimal cost function in the vector of installed storage capacities, which reveals that the marginal value of storage capacity is largest for initial investments. Then, we provide an explicit characterization of the sensitivity of the optimal cost to an initial investment in storage capacity at a node in the network -- which we refer to as the locational marginal value of storage capacity. Our main result (Theorem \ref{thm:mv}) offers an expression for this marginal value in terms of nodal price expectations. Computing these price expectations is an area of active research. See \cite{ji2015probabilisticCONF, ji2015probabilisticJOURNAL} for example. In essence, our result provides a tool to compute the marginal value of initial investment in storage capacity in the power grid. The marginal value, in turn, provides a first-order approximation to the value of storage capacity that is challenging to rigorously characterize. To provide further insights, we also offer an upper bound on the marginal value in terms of the total variation of the same price process, and identify specific network and cost structures in  Proposition \ref{prop:spCase}, where this upper bound is achieved.

The paper is organized as follows. We begin by describing the system model in Section \ref{sec:models} and formulate the SO's problem in Section \ref{sec:multiED}. Then, in Section \ref{sec:lmp}, we introduce the concept of nodal pricing, and establish certain parametric properties of nodal prices that will prove central to establishing our main results in Section \ref{sec:results}. We conclude the paper with a  detailed analysis of a two-node power network in Section \ref{sec:2Node} to further interpret and illustrate our results. All proofs are contained in the appendix.

\emph{Notation:} Let $\Rset$ (resp. $\Rset_+$) denote the set of real (resp. nonnegative) numbers. For $x \in \Rset$, let $x^+ := \max\{x, 0\}$. 
For a vector $x$, let $x^i$ denote its $i$-th entry and $x^\top$ denote its transpose. Similarly, for a matrix $X$, let $X^{ij}$ denote the entry at the $i$-th row and the $j$-th column of $X$, and  $X^\top$ denote its transpose.
For a sequence of elements $x:= (x_0, x_1, \ldots)$, define $x_{\leq k} := (x_0, \ldots, x_k)$. 
For $h: \Rset^n \to \Rset$, let $\partial h/\partial x^i \vert_{x=x_0^+} $ denote the right-hand derivative of $h$ with respect to the $i$-th coordinate at $x_0$. 
For a random variable $X$, let $\E[X]$ denote its expectation. 
For an event $\mathcal{E}$, define $\prob\{\mathcal{E}\}$ as the probability of that event. 
For two sets $\Acal$ and $\Bcal$, define $\Acal \cap \Bcal$, $\Acal \cup \Bcal$, and $\Acal \setminus \Bcal$ as the intersection, union, and the set difference of $\Acal$ and $\Bcal$, respectively. 
For a set $\Acal$ in Euclidean space, define $\interior{\Acal}$ as its interior and $\bound{\Acal} = \Acal \setminus \interior{\Acal}$ as its boundary.


\section{System Model}
\label{sec:models}

In Definition \ref{def:LMVstorage}, we define the locational marginal value of storage capacity as the sensitivity of the optimal value of the multi-period economic dispatch (ED) problem to an increase in storage capacity at a particular node in the network. We begin by describing the basic system components required to formulate the multi-period economic dispatch problem.

\subsection{Network Model} 
\label{sec:network}
Consider a power network described by a connected undirected graph on $m$ nodes (or buses), labeled $1,2, \ldots, m$, and $\ell$ edges (or transmission lines). The set of feasible nodal power injections is defined as 
\begin{align}
\label{eq:injPoly}
\Pcal  &:=  \{ x \in \Rset^m \ \vert \ x = Y \theta,  -f \leq B \theta \leq f \text{ for some }\theta \in \Rset^{m}\},
\end{align}
We refer to the set $\Pcal$ as the \emph{injection polytope}. 
Here, $Y \in \Rset^{m \times m}$ denotes the bus admittance matrix, defined by $Y^{ij} := -y^{ij}$ for $i \neq j$, and $Y^{ii} := y^{ii} + \sum_{j \neq i} y^{ij}$, where $y^{ij}$ is the susceptance of the transmission line joining buses $i$ and $j$, and $y^{ii}$ is the susceptance of the shunt element at bus $i$. The matrix
 $B \in \Rset^{\ell \times m}$ is the (weighted) incidence matrix of the network.
 For a transmission line $k$ joining buses  $i$ and $j$, define 
  $B^{ki} = -B^{k j} := y^{i j}$ and $B^{kr} = 0$ for all $r \neq i$ or $j$. Finally, we denote by $f \in \Rset^\ell_+$  the vector of transmission line power capacities, and  $\theta \in \Rset^m$  the vector of bus voltage phase angles.
The above description of the set of feasible injections is derived using a linear approximation of the Kirchhoff's laws, 
commonly known as the \emph{DC power flow model}. This widely used linear model assumes
purely reactive transmission lines, constant bus voltage magnitudes, and small voltage phase angle differences across transmission lines. See \cite[Ch. 6]{andersson2004modelling}, \cite[Ch. 9]{grainger1994power} for a detailed derivation.


\subsection{Net Demand Process}
\label{sec:netDemProc}
We assume time to be discrete 
and consider the operation of the power network over a finite horizon of $N$ time periods indexed by $k=0,1, \dots, N-1$. Recall that by net demand, we mean the nominal (inelastic) demand less any variable supply. We denote the vector of nodal net demands across the network at time $k$ by  $\xi_k := (\xi_k^1, \xi_k^2, \ldots, \xi_k^m)^{\top} \in \Rset^m$. 
The system operator (SO) aims to balance the net demand using dispatchable resources at minimum cost. The sign convention is such that $\xi_k^i \leq 0$ represents a net supply of energy, while $\xi_k^i > 0$ indicates a net consumption of energy at node $i$ in period $k$. The spatio-temporal evolution of net demand is modeled as a discrete time vector random process $\{\xi_0, \xi_1, \dots,\xi_{N-1} \}$, where $\xi_k$ takes values in the set $\Xi \subset \Rset^m$ for each $k = 0, \ldots, N-1$. Assume that the joint distribution of this random process is known. It is important to note that we do not require the process to be stationary or  independent across time or space.

\subsection{Cost Structure} 
\label{sec:genCost}
Our model is such that each bus in the power network is allowed to have both dispatchable generation and elastic demand.
This is reflected in the \emph{nodal cost of generation} function $g^i : \RR \rightarrow \RR$ for $i = 1, \dots, m$. 
Specifically, the cost of producing $v^i$ amount of power at bus $i$ is defined as
$$ g^i(v^i) := \alpha^i (v^i)^+ - \beta^i(-v^i)^+.$$
When $v^i$ is positive, $g^i(v^i)$ denotes the cost of generating $v^i$. When $v^i$ is negative, $-g^i(v^i)$ denotes the utility of consuming $-v^i$.  The quantities $\alpha^i$ and $\beta^i$ represent the marginal cost of generation and the marginal utility of consumption at node $i$, respectively. 
We assume that $\alpha^i \geq \beta^i \geq 0$ for all $i$. Such cost structure mirrors the formulation in \cite{wu1996folk}.
It follows that the nodal cost function $g^i$ is convex, nondecreasing,  piecewise linear, and independent of time. The results stated in this paper are easily generalized to the case in which $g^i$ is allowed to vary with time and has an arbitrary, but finite number of break-points.
Henceforth, we will refer to the  vector $v = (v^1,\dots, v^m)^\top$ as a \emph{dispatch}. The system-wide cost incurred by a dispatch $v \in \Rset^m$ is thus given by
$$ g(v) := \sum_{i=1}^m g^i (v^i).$$

\subsection{Energy Storage Model} 
\label{sec:storage}

Consider a collection of $m$  {perfectly efficient} energy storage devices built into the network,  where we allow at most a single storage device at each node $i=1,\dots,m$. The collective storage dynamics can be modeled according to the following linear difference equation
\beq
\label{eq:storageDynamics} z_{k+1} = z_k  + u_k, \quad k=0,\dots,N-1,
\eeq
where $z_k \in \Rset_+^m$ denotes the vector of energy storage states, just preceding period $k$, and $u_k \in \Rset^m$  denotes the vector of energy extractions or injections during period $k$. We adopt the convention, where $u_k^i \geq 0$ (resp. $u_k^i < 0$) represents a net energy injection into (resp. extraction from) the storage device at node $i$ during period $k$. Without loss of generality, assume a zero initial condition ($z_0 = 0$) throughout. 
For a vector of energy storage capacities $b \in \Rset_+^m$ installed across the network, the storage dynamics are constrained as
\begin{align*}
z_k \in \Zcal(b), \text{ where } \Zcal(b) := \{ z \in \Rset_+^m \ \vert \ 0 \leq z \leq  b \}  
\end{align*}
for all $k=0,\ldots,N$. 
This work ignores ramping constraints on the incremental injections/extractions or round-trip inefficiency/dissipative losses in the storage devices. We refer the reader to Remark \ref{rem:nonIdeal} for a discussion on incorporating such non-idealities into the analysis.


\section{Problem Formulation}

\label{sec:multiED}
In what follows, we formalize the SO's problem of multi-period economic dispatch (ED) with storage as a finite-horizon, constrained  stochastic control problem with perfect state feedback.
Working within this setting, we define, in Section \ref{sec:lmv}, the locational marginal value of storage capacity as the parametric sensitivity of the optimal value of the multi-period ED problem.

\subsection{Multi-period Economic Dispatch Problem}
\label{sec:multiED.1}
For each time period $k = 0, \ldots, N-1$, define $(z_k, \xi_k) \in \Rset^m_+ \times \Rset^m$ as the system state. Recall that $z_k$ is the vector of energy storage states just preceding period $k$, and $\xi_k$ is the vector of net demands across the network at period $k$. 
It is natural to assume that the SO has perfect state feedback. Thus, the available information at time $k$ is defined as 
$I_k := (z_{\leq k}, \xi_{\leq k})$.
It is the SO's task to determine a \emph{control policy} $\vpi := \left( (\mu_0, \nu_0), \ldots, (\mu_{N-1}, \nu_{N-1}) \right)$, which maps  available information $I_k$  to inputs $(u_k, v_k)$ at each time $k$. Namely,  $u_k = \mu_k (I_k)$ and $v_k = \nu_k (I_k)$ for each $k = 0, \ldots, N-1$. 
To ensure that a control policy respects the network and storage capacity constraints at each time, 
we define the following notion of admissibility. 

\begin{definition}[Admissible policies] Given a vector of storage capacities $b \in \Rset^m_+$,
a control policy $\v{\pi} =\left( (\mu_0, \nu_0), \ldots, (\mu_{N-1}, \nu_{N-1}) \right)$ is said to be \emph{admissible}, if 
$$z_k + \mu_k(I_k) \in \Zcal(b) \quad and  \quad \nu_k(I_k) - \mu_k(I_k) - \xi_k \in \Pcal,$$
almost surely for each $k=0,\dots,N-1$. 
Denote by $\Pi(b)$, the space of all admissible policies.
\end{definition}

The {expected cost of dispatch} over $N$ periods under an admissible control policy $\vpi \in \Pi(b)$ is then given by
\begin{align}
\label{eq:cost}
J^{\v{\pi}}(b) \ := \ \E \left[  \ \ \sum_{k=0}^{N-1} g( \nu_k(I_k) ) \ \right],
\end{align}
where the expectation is computed with respect to the known distribution on the net demand process. The SO seeks an admissible control policy that minimizes the above cost.

\begin{definition}[Optimality criterion] Given a vector of storage capacities $b \in \Rset^m_+$, the minimum expected total cost of dispatch over $N$ periods is defined as $ J^*(b) := \text{infimum} \left\{ J^{\v{\pi}}(b) \ \vert \ {\v{\pi} \in \Pi(b)} \right\}$.
An admissible control policy $\v{\pi}^* \in \Pi(b)$ is said to be \emph{optimal}, if $J^{\v{\pi}^*}(b) = J^{*}(b)$.
\end{definition}

\subsection{Defining the Value of Storage}
\label{sec:lmv}

The capital investment cost of a grid scale storage asset depends heavily on its energy capacity \cite{eyer2010energy}. Thus, in order to correctly size and place such assets within a given power network, it is critical to quantify the maximum benefit one might derive from their utilization in multi-period ED. A natural way to measure the value of storage capacity within the context of our model is the maximum reduction in dispatch cost achievable with a collection of storage assets described by $b \in \Rset^m_+$ -- in other words, $J^*(0) - J^*(b)$. 
We have the following structural result on $J^*(b)$.
\begin{theorem}
\label{thm:convex} For each $b\in \Rset_+^m$, there exists an optimal control policy $\v{\pi}^* \in \Pi(b)$. Furthermore, $J^*(b)$ is convex and non-increasing in $b$.
\end{theorem}
\rev{We remark that the proof follows from standard arguments on the existence of an optimal control policy in a stochastic control problem and on Jensen's inequality. It is omitted for brevity.}

Theorem \ref{thm:convex} reveals that $J^*(b)$ is convex, which in turn implies  diminishing returns on investment in storage capacity. As a result, the greatest \emph{marginal value of storage capacity} is derived from \emph{initial investments}.  \rev{Adoption of energy storage technology in the power grid is currently in its infancy. The value of initial investment in storage capacity at various nodes in a power network will then serve to inform the storage adoption decisions for system operators in practice. Calculating this value for an arbitrary storage capacity $b \in \Rset^m_+$ hinges on characterizing the optimal storage control policy $\pi^*(b)$. Theorem \ref{thm:convex} ensures the existence of such a policy; its explicit characterization, however, remains an open question.}\footnote{\rev{It is well-known that optimal control policies are difficult to characterize in stochastic control problems that include constraints that involve both the states and the inputs.}} \rev{The locational marginal value of initial investment in storage capacity, as defined next, provides a first-order approximation of $J^*(b)$ near the origin. 
Succinctly put, the nodes with larger marginal values at the origin represent better choices for initial siting of storage assets. Our main result in Section \ref{sec:results} offers a tool to estimate these marginal values empirically.}

\begin{definition}[Locational marginal value] \label{def:LMVstorage} For each $i = 1,\ldots,m$, the \emph{locational marginal value} (LMV) of initial investment in storage capacity at node $i$ is defined as 
\begin{align*}
{\sf LMV}^i = -\frac{\partial J^*(b)}{\partial b^i}\bigg\vert_{b =0^+} .
\end{align*}
\end{definition}
The convexity of $J^*(b)$ in $b \in \Rset_+^m$ (from Theorem \ref{thm:convex}) ensures that the coordinate-wise right-hand partial derivatives in the above definition exist.


\section{Locational Marginal Pricing}
\label{sec:lmp}

We now define a stochastic price process that will prove essential in characterizing the locational marginal value of storage capacity in Theorem \ref{thm:mv}. This price process, being endogenously defined,  is linked to the net demand process through the dual optimal solution of a multi-parametric linear program, introduced next. We first require the following definition of single-period economic dispatch (ED).

\subsection{Single-period Economic Dispatch problem}
Given a vector of net demands $\xi \in \Rset^m$, define the single-period ED problem and its optimal cost as
\begin{align}\label{eq:ED}
Q(\xi) := \underset{v \in \Rset^m}{\text{infimum  }} g(v), \ \text{subject to} \ \ v - \xi \in \Pcal.
\end{align}
The above multi-parametric linear program is feasible and the infimum is achieved for each $\xi \in \Rset^m$. It follows from the definition of the polytope $\Pcal$ in \eqref{eq:injPoly} that the constraint $v-\xi \in \Pcal$ can be equivalently represented as
 $v - \xi = Y \theta$ and $-f \leq B\theta \leq f$ for some  $\theta \in \Rset^m$. Let $\lambda (\xi) \in \Rset^m$ be the optimal Lagrange multiplier associated with the power balance constraint $v - \xi = Y\theta$. We have parameterized the Lagrange multiplier by $\xi \in \Rset^m$ to make explicit its dependence on the vector of net demands.

\rev{Recall that $\Xi$ denotes the support of the net demands at each period. Important to the sequel are structural properties of the parametric optimal value $Q(\xi)$, and the parametric dual optimal solution $\lambda(\xi)$ over elements in $\Xi$ and their neighborhoods. To facilitate such analyses, we characterize $Q(\xi), \lambda(\xi)$ in Lemma \ref{lemma:QxiLxi} over a full-dimensional polytope $\Kcal$ containing $\Xi$ in its interior (denoted by $\interior{\Kcal}$).}\footnote{\rev{We say $\Kcal$ is full-dimensional to indicate that $\interior{\Kcal}$ is non-empty.}} \rev{The choice of such sets is arbitrary. Henceforth, we fix $\Kcal$ for the ease of exposition. Lemma \ref{lemma:QxiLxi} follows largely from arguments in \cite[Theorem 7.2]{BorrelliBook}.
Its proof is deferred till Appendix \ref{app:proof:QxiLxi}.}

Stating Lemma \ref{lemma:QxiLxi} requires an additional notation.  For a polyhedral set $\Kcal$, a finite collection of sets $\Lcal$ defines a \emph{polyhedral partition} of $\Kcal$, if (i) each set  $\Scal \in \Lcal$ is polyhedral, (ii) $\bigcup_{\Scal \in \Lcal} \Scal = \Kcal$, and (iii) $\interior{\Scal} \cap \interior{\Scal^{'}}$ is empty for two distinct sets $\Scal, \Scal^{'}$ in $\Lcal$.

\begin{lemma}
\label{lemma:QxiLxi}
\rev{Let $\Kcal  \subset \Rset^m$ be a full-dimensional polytope such that $\Xi \subset \interior{\Kcal}$. Then, there exists a polyhedral partition $\Lcal = \{\Scal_1, \dots,\Scal_{|\Lcal|}\}$ of $\Kcal$ that satisfies:}
\begin{enumerate}[(i)] \setlength{\itemsep}{.08in}
\item $Q(\xi)$ is affine over $\Scal_\ell$ for each $\ell = 1,\dots, |\Lcal|$. 
\item $\lambda(\xi)$ is nonnegative and constant and $\nabla Q(\xi) = \lambda(\xi)$ over  $\interior{\Scal_\ell}$ for each $\ell = 1,\dots, |\Lcal|$.
\item The union of the boundaries of the sets comprising the polyhedral partition, i.e., 
$\Bcal := \bigcup_{\ell = 1}^{|\Lcal|} \partial \Scal_\ell$,
has zero Lebesgue measure.
\end{enumerate}
\end{lemma}

We denote the polyhedral, piecewise constant dual multiplier specified in Lemma \ref{lemma:QxiLxi} as
\begin{align} \label{eq:pwafunc}
\lambda(\xi) = p_{\ell}, \quad \text{if} \ \ \xi \in \interior{\Scal_\ell},
\end{align}
where $p_\ell \in \Rset^m$ for all $\ell =1, \dots, |\Lcal|$.


\begin{figure}[t]
\centering
\includegraphics[width=0.25\textwidth]{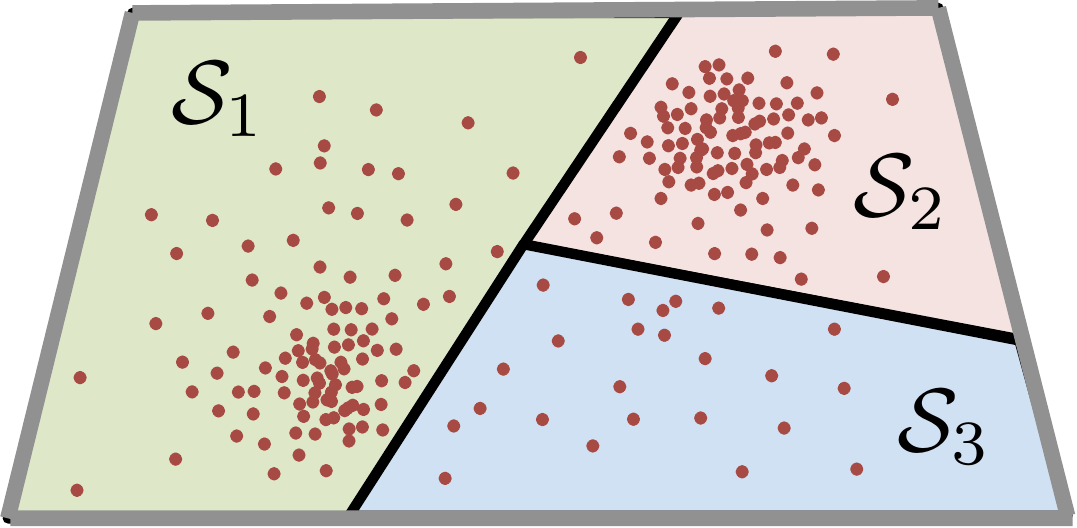} \label{fig:polyPart} 
\caption{An example with $m=2$. The collection of dots represents $\Xi$. The polytope with gray boundary represents $\Kcal$. Sets $\Scal_1, \Scal_2$, and $\Scal_3$ define the polyhedral partition $\Lcal$ of $\Kcal$. The collection of black and gray line segments represent $\Bcal$. 
} 
\label{fig:polyPart}
\end{figure}

\subsection{The Locational Marginal Price Process}

\rev{Lemma \ref{lemma:QxiLxi} implies that $\lambda(\xi)$ is uniquely defined over ${\Kcal}$, except possibly on the union of the boundaries of its polyhedral partition, denoted by $\Bcal$. We make the following technical assumption in order to ensure that $\lambda(\xi)$ is well-defined on all elements in $\Xi$.
\begin{assumption}
\label{ass:noS} The set $\Xi$ is finite and $\Xi \cap \Bcal$ is empty.\footnote{\rev{On first glance, it may appear that the statement in Assumption \ref{ass:noS} depends on the particular choice of $\Kcal$ containing $\Xi$ in its interior. It can be shown, however, that $\lambda(\xi)$ is uniquely defined over $\Rset^m$, except on a set of zero Lebesgue measure. As long as $\Xi$ does not intersect this zero measure set, any $\Kcal$ containing the finite set $\Xi$ in its interior will satisfy that $\Xi \cap \Bcal$ is empty.}}
\end{assumption}
Figure \ref{fig:polyPart} illustrates the sets $\Kcal$, $\Lcal$, $\Bcal$ and $\Xi$ for an example that satisfies Assumption \ref{ass:noS}.} \rev{Utilizing the parametric dual optimal solution of \eqref{eq:ED}, we next define a nodal price process in terms of which we characterize the locational marginal value of storage in the next section.}
\begin{definition}[Nodal price process] \label{def:nodalprice}
Suppose Assumption \ref{ass:noS} holds. For $i=1, \ldots, m$, the  \emph{price process} at node $i$ is defined as the scalar random process ${\lambda}^i := \{ \lambda^i_0, \ldots, \lambda^i_{N-1} \}$, where 
\begin{align} \label{eq:nodprice}
\lambda^i_k := \lambda^i (\xi_k), \ \ k=0\ldots,{N-1},
\end{align}
and $\lambda(\cdot)$ is the optimal Lagrange multiplier associated with the power balance constraint in problem \eqref{eq:ED}.
\end{definition}

The definition of nodal price according to $\lambda_k^i = \lambda^i (\xi_k)$ coincides  with the standard approach of nodal or locational marginal pricing in wholesale electricity markets  today \cite{schweppe1988spot, hogan1992contract}. Specifically, the nodal price $\lambda_k^i$ equals the marginal cost of serving an additional unit of demand at node $i$ and time $k$, in the absence of energy storage. The vector of nodal prices at time $k$ is denoted by $\lambda_k := \left( \lambda_k^1, \ldots, \lambda_k^m \right)^\top$.

The nodal price process is clearly stochastic, given its explicit dependency on the underlying net demand process. Its one-step look-ahead predictor will prove important in the sequel. We define it as
\begin{align} 
{\lambda}_{k+1 | k} := \E [\lambda_{k+1} \  \vert \  \xi_{\leq k}] \in \Rset^m, \label{eq:condMean}
\end{align}
for $k=0, \ldots, N-2$. 
At time $k$, the predictor ${\lambda}_{k+1 | k}$  denotes the expected value of the nodal prices at the following time $k+1$, conditioned on the history of the net demand through time $k$.


\section{Locational Marginal Value of Storage}
\label{sec:results}

Theorem \ref{thm:mv} contains our main result,
which  offers an explicit characterization of the locational marginal value of energy storage in terms of a certain measure of nodal price variation. We first require a definition. 
The \emph{total variation} of a scalar sequence  $x = (x_0, \dots, x_{N-1})$ is defined as
\begin{align}
\label{eq:TV}
\mathsf{TV}(x) \ := \ \sum_{k=0}^{N-2} \left\vert x_{k+1} - x_k\right\vert.
\end{align}

\begin{theorem}
\label{thm:mv}
Suppose Assumption \ref{ass:noS} holds. For $i=1,\ldots,m$, the marginal value of initial investment in storage capacity at node $i$ is given by
\begin{align}
{\sf LMV}^i = \E \left[ \sum_{k=0}^{N-2}({\lambda}^i_{k+1 | k} - \lambda^i_k)^+ \right].
\label{eq:mv}
\end{align}
Moreover, it is bounded from above as
\begin{align}
%
{\sf LMV}^i   \leq  \frac{1}{2} \ \E \left[ \mathsf{TV}(\v{\lambda}^i) \right]  +  \frac{1}{2} \ \E \left[ \lambda^i_{N-1} - \lambda^i_{0} \right]. 
\label{eq:mvUp}
\end{align}
\end{theorem}
We defer the proof of Theorem \ref{thm:mv} to Appendix \ref{app:proof:mv}. The formulae offered by Theorem \ref{thm:mv} admit intuitive \emph{dual interpretations} as the maximum expected revenues achievable through price arbitrage with storage. More specifically, consider a setting in which a storage owner-operator seeks to dispatch her storage device with the objective of maximizing her expected revenue through arbitrage against the sequence of stochastic nodal prices $\{\lambda^i_0, \ldots, \lambda^i_{N-1}\}$. Assuming the nodal price process to be unaffected by the storage owner-operator's control actions, one can show the optimal causal arbitrage policy to be of \emph{price threshold-type}. That is, at each time period $k$, one compares the current price $\lambda_k^i$ with a threshold given by the one-step look-ahead expected  price $\lambda_{k+1|k}^i$ \rev{as defined in \eqref{eq:condMean}} . If  the price is expected to increase (i.e., $\lambda_{k+1|k}^i \geq  \lambda_{k}^i$), the optimal policy dictates that one buys an amount of energy that fills the storage device to capacity. If, on the other hand, the price is expected to fall (i.e., $\lambda_{k+1|k}^i <  \lambda_{k}^i$), it is optimal to  sell an amount of energy that fully empties the storage device. And naturally, it is always optimal to empty the storage device at the terminal stage. It follows that, given a storage capacity of $b \geq 0$, the expected revenue achieved under such a policy is equal to $b\cdot {\sf LMV}^i$, as defined  by equation \eqref{eq:mv}. And, one can show that the upper bound in \eqref{eq:mvUp} is similarly derived from the maximum expected revenue achievable with perfect foresight of the nodal price process. 

Ultimately, Theorem \ref{thm:mv} reveals the value of initial investment in storage capacity at a particular location in the power network to depend on the variation in the net-demand process insofar as it manifests itself as variation in the corresponding nodal price process\footnote{The former does not necessarily imply the latter. This is made clear by the example studied in Section \ref{sec:2Node}.}  -- a polyhedral, piecewise constant function of net-demand (cf. Lemma \ref{lemma:QxiLxi}).
We refer the reader to the parametric analysis of a two-node power network in 
Section \ref{sec:2Node}, which illustrates the effect of network transmission capacity on the behavior of this mapping. 

\rev{If the cost functions $g^i$'s are smooth, one can obtain piecewise linear approximations to such $g^i$'s with arbitrary precision. Our results then provide a tool to approximate the locational marginal value of storage. Guarantees on the approximation quality for the marginal value can be obtained in terms of the accuracy of the piecewise linear approximations to $g^i$'s.}

\begin{remark}[Calculating the locational marginal value]
\label{rem:calcMarginal}
Theorem \ref{thm:mv} shows that  the calculation of the locational marginal value of storage reduces to the calculation of nodal price expectations. Of particular importance is the one-step look-ahead conditional expectation of nodal prices given by $\lambda_{k+1|k}  =  \E [\lambda(\xi_{k+1}) \  \vert \  \xi_{\leq k}] $.  Using the  polyhedral, piecewise constant representation of  $\lambda(\cdot)$ specified in \eqref{eq:pwafunc}, we arrive at the following simplified form
\begin{align*}
\lambda_{k+1|k}   =   \sum_{\ell=1}^{|\Lcal|}  \  p_{\ell} \cdot \prob\{\xi_{k+1} \in \Scal_\ell \ | \ \xi_{\leq k} \} ,
\end{align*}
where recall 
that $p_\ell \in \Rset^m$ denotes the nodal price vector induced by a net-demand vector belonging to the polyhedral set $\Scal_\ell$. 
The challenge in computing  $\lambda_{k+1|k}$ thus reduces to the calculation of the conditional probabilities $\prob\{\xi_{k+1} \in \Scal_\ell \ | \ \xi_{\leq k} \}$. We refer the reader to recent work \cite{ji2015probabilisticCONF, ji2015probabilisticJOURNAL} that is dedicated precisely to the resolution of this challenge, and offers a detailed exposition into the analytical and empirical calculation of such conditional probabilities. 

\end{remark}

\begin{remark}[Modeling non-idealities in storage]
\label{rem:nonIdeal} 
Theorem \ref{thm:mv} assumes perfectly efficient storage assets with no ramping constraints. Storage devices in practice, however, suffer from non-idealities like dissipative losses, roundtrip efficiency losses, and have limited ramping capabilities. With each of these non-idealities, one can show that $b \cdot {\sf LMV}^i$ still equates to the maximum expected revenue a storage owner-operator can derive from a causal arbitrage against the price process at node $i$ with a storage device of capacity $b \geq 0$. For a dissipative storage model, said maximum expected revenue equals $\E \left[ \sum_{k=0}^{N-2}(\gamma {\lambda}^i_{k+1 | k} - \lambda^i_k)^+ \right]$, where $z_{k+1} = \gamma z_k + u_k$ at time $k$ for a dissipation rate $\gamma \in (0,1)$. 
However, \rev{when a storage device has a roundtrip efficiency loss or  bounded ramp rates, said maximum expected revenue does not admit a succinct representation}.
\end{remark}

\subsection{Achieving the Upper Bound}
\label{sec:upBound}

Theorem \ref{thm:mv} offers formulae to enable the tractable calculation and upper bounding of the locational marginal value of storage, ${\sf LMV}^i$. In the following result (Proposition \ref{prop:spCase}), we identify sufficient conditions under which our upper bound is achieved. Its proof can be found in Appendix \ref{app:proof:spCase}.
\begin{proposition}
\label{prop:spCase}
Suppose that Assumption \ref{ass:noS}  and the following conditions hold: 
\begin{enumerate}[(a)]
\item the graph of the power network is acyclic, and
\item the energy costs are spatially homogeneous, i.e., $\alpha^j = \alpha$ and $\beta^j = \beta$ for all $j=1,\dots,m$ for some $\alpha \geq \beta \geq 0$.
\end{enumerate}
Then, for each $i = 1,\ldots, m$, we have that
\begin{align*}
\lambda^i(\xi) \in \{ \alpha, \beta\} \ \ \text{for all} \ \xi \in \Xi,
\end{align*}
and
\begin{align*}
\mathsf{LMV}^i &= \frac{1}{2} \E \left[ \mathsf{TV}(\v{\lambda}^i) \right]  +  \frac{1}{2}  \E \left[ \lambda^i_{N-1} - \lambda^i_{0} \right],
\end{align*}
Further, $\mathsf{LMV}^i$ equals $(\alpha - \beta)$ times the expected number of periods $k$ for which $\lambda^i_{k}=\beta$ and $\lambda^i_{k+1}=\alpha$.
\end{proposition}
We shed light on the meaning of Proposition \ref{prop:spCase}. 
Consider again the setting in which an arbitrageur seeks to operate a storage device located at node $i$ to maximize the expected revenue she derives through the buying and selling of energy against the nodal prices $\lambda^i_0, \ldots, \lambda^i_{N-1}$. Recall that -- assuming the nodal prices to be unaffected by the actions of the storage owner-operator -- the maximum expected revenue achievable with a storage device of capacity $b \geq 0$ is  given by  $b \cdot {\sf LMV}^i$.
With perfect foresight into said prices, she can garner $b$ times the upper bound in \eqref{eq:mvUp}. To appreciate where the gap between the revenues stems from, consider the following control policy the storage owner-operator implements for optimal price arbitrage under perfect foresight. Buy energy to charge the device to capacity, whenever the nodal price process is at its local minimum. Then, sell to empty the device at the following local maximum. Also, always sell to empty it at the last period. Notice that such a policy, in general, cannot be executed causally. Deciding whether a scalar stochastic process is currently at a local extremum requires foresight into its future, in general. We, however, circumvent this difficulty when the network is acyclic and the costs are spatially homogeneous, as $\lambda^i_k \in \{\alpha, \beta \}$ for all $i=1,\ldots,m$ and $k=0,\ldots,N-1$ under the assumptions of Proposition \ref{prop:spCase}.
That is, the nodal prices can only be either $\alpha$ or $\beta$.\footnote{\rev{Nodal prices can take values other than $\alpha$ or $\beta$, if the network contains cycles.}} Consequently, the nodal price process is at a local minimum, whenever the current nodal price is $\beta$, and at a local maximum, whenever it is $\alpha$. As a result, one can causally implement the optimal control policy with perfect foresight. In turn, ${\sf LMV}^i$ achieves its upper bound.


\section{Analysis of a Two-Node Network}
\label{sec:2Node}
\begin{figure*}[htb]
\centering
\subfigure[$\lambda^1(\xi)$: price at node 1 as a function of $\xi$. ]{ { \scalebox{0.5}{\includegraphics*{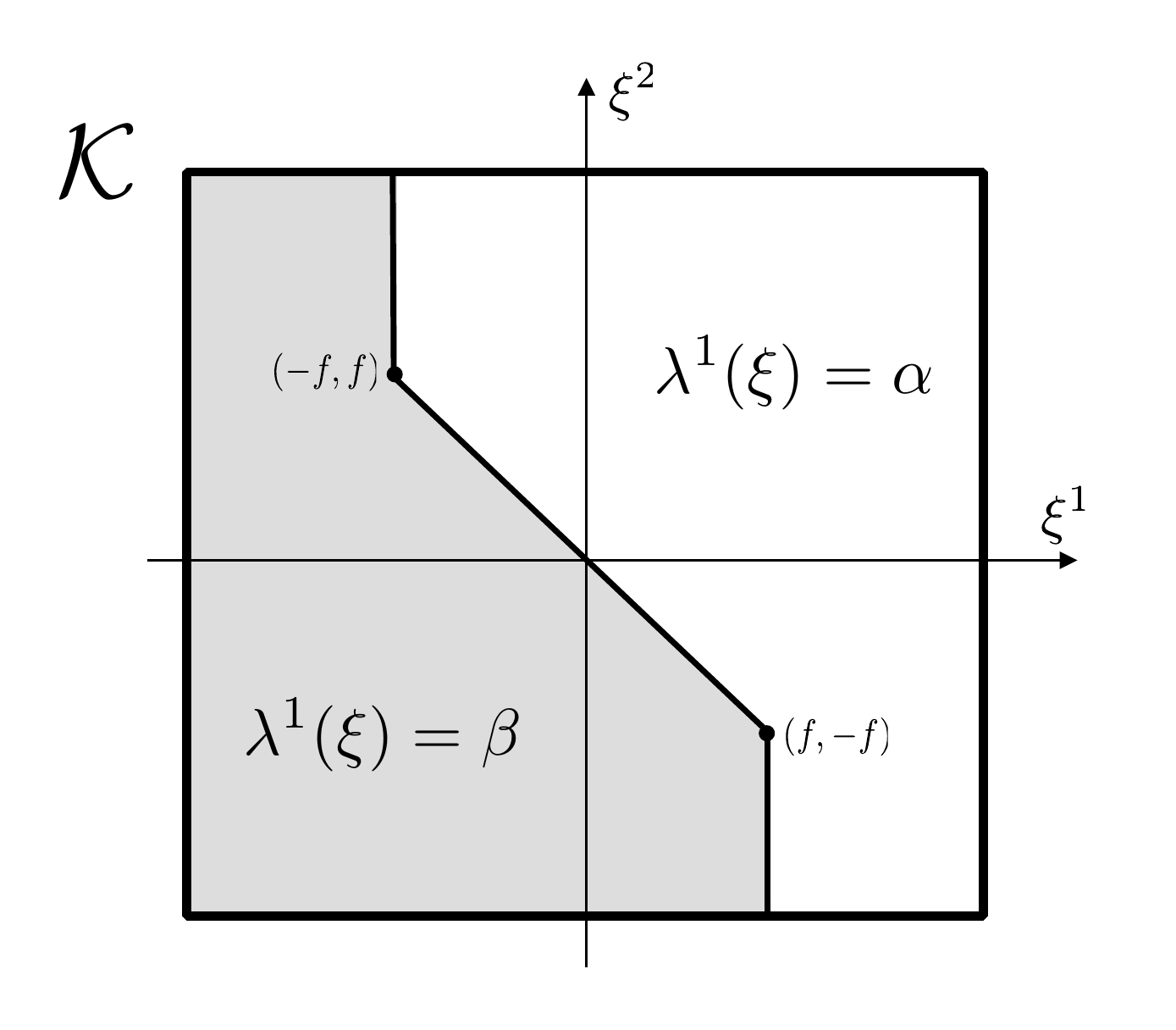}}} \label{fig:2Node.1} } \ \ \ \ \ \ \ 
\subfigure[$\lambda^2(\xi)$: price at node 2 as a function of $\xi$.]{ {\scalebox{0.5}{\includegraphics*{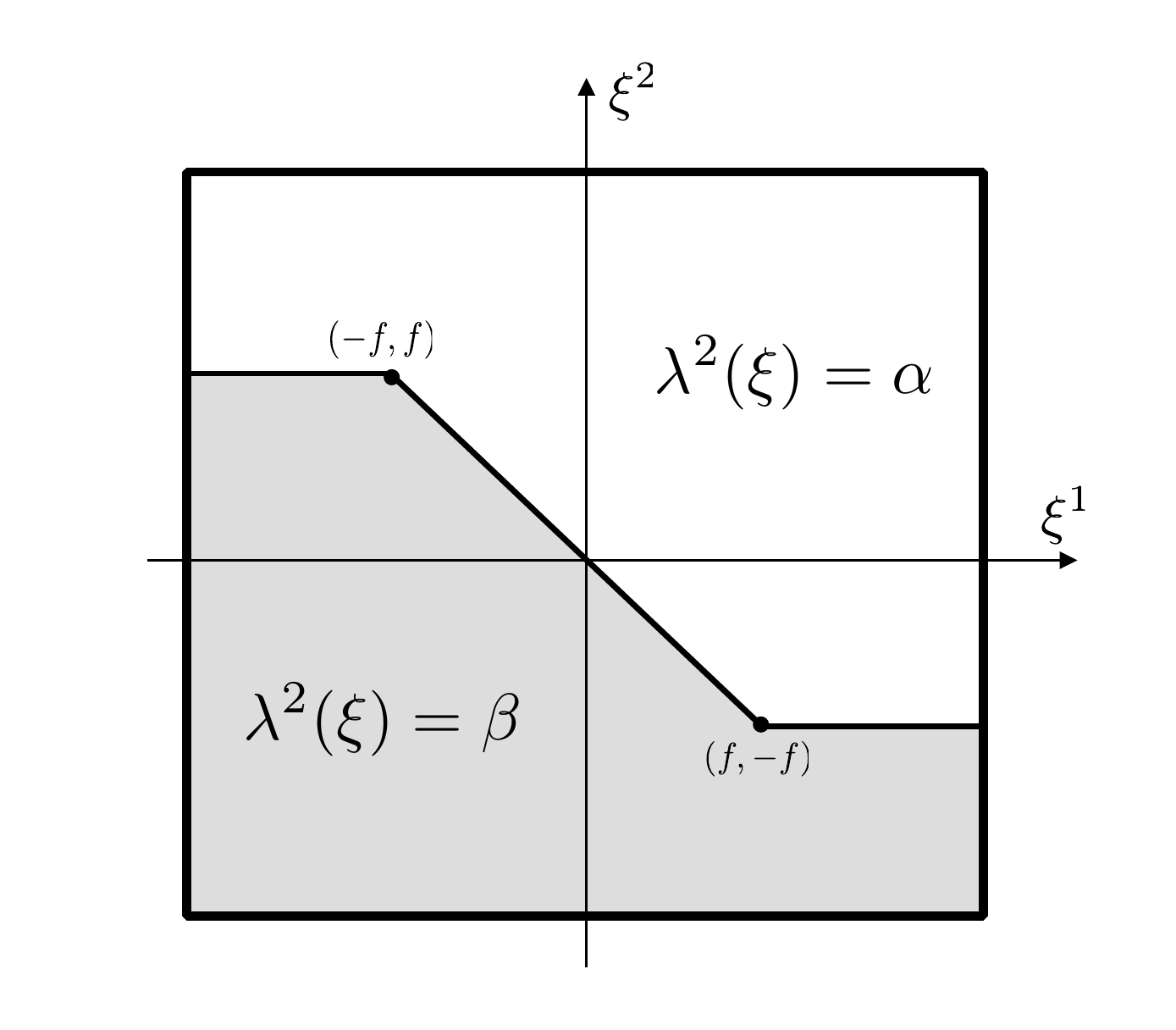}}} \label{fig:2Node.2}}
\caption{This diagram depicts the nodal prices (i.e., optimal Lagrange multipliers) as parametric functions of the net demand $\xi := (\xi^1, \xi^2)$ over $\mathcal{K}$ (which is taken to be a closed square centered at the origin). For each node $i = 1,2$, it holds that $\lambda^i(\xi) = \beta$ for net demand  values  in the shaded region,  and $\lambda^i(\xi) = \alpha$ for net demand  values  in the white region.} \label{fig:2Node}
\end{figure*}

\rev{
We now analyze the special case of a two-node network to illuminate the effect which the network transmission capacity has upon the locational marginal value (LMV) of storage, as revealed by our theoretical results. More precisely, consider a network with two nodes (labeled 1 and 2) joined by a single transmission line having capacity $f \geq 0$. And, let conditions (a)-(b) of Proposition \ref{prop:spCase} hold. It follows that the nodal prices can take one of two values, $\alpha$ or $\beta$, depending on the value of net demand $\xi$. We depict this parametric dependency of nodal prices on net demand in Figure  \ref{fig:2Node}, which indicates that  $\lambda^i(\xi) = \beta$ for net demand  values  in the shaded region, and  $\lambda^i(\xi) = \alpha$ for net demand  values  in the white region, for each node $i=1,2$.
}  

\rev{ 
Figure \ref{fig:2Node} also reveals a precise relationship between the LMV of storage capacity and the transmission line capacity.  For example, upon examination of both Figures  \ref{fig:2Node.1}-\ref{fig:2Node.2}, it becomes  immediate to see that the LMV of storage capacity at each node is directly proportional to the expected number of times at which the net demand process   $(\xi_0, \dots, \xi_{N-1})$ transitions from the shaded (low price) region to the white (high price) region.  And,  as can be directly inferred from the figures, the frequency with which the net demand process exhibits such crossings depends explicitly on the transmission line capacity, and the extent to which it promotes the `mixing' of net demands between the two nodes.  We shed light on this dependency by considering the limiting cases of \emph{low} ($f \rightarrow 0$) and \emph{high} ($f \rightarrow \infty$) transmission capacity. }

 \rev{
First, as one might naturally expect, a reduction in the transmission capacity between nodes 1 and 2 serves to attenuate the 
 degree to which fluctuations in net demand at node 1 have an influence on the marginal value of storage at node 2, and vice versa. In the limit as $f \rightarrow 0$, a straightforward calculation reveals the LMV at each node $i$ to simplify to 
 \begin{align*}
\lim_{f \rightarrow 0} {\sf LMV}^i  =   (\alpha - \beta) \sum_{k=1}^{N-1} \mathbb{P}\{ \xi^i_{k-1} < 0, \ \xi^i_k > 0 \}.
\end{align*}
In words,  the marginal value of storage at each node $i$ becomes dependent only on the statistical variation of its local net demand process $\{\xi_k^i\}$, as measured through its  \emph{expected number of zero-upcrossings}. 
}

\rev{
In the limit as the transmission capacity grows large, an analogous argument reveals the marginal value of storage at each node $i$ to satisfy
\begin{align*}
\lim_{f \rightarrow \infty} {\sf LMV}^i  =   (\alpha - \beta) \sum_{k=1}^{N-1} \mathbb{P}\{ \xi^1_{k-1} + \xi^2_{k-1}  < 0, \ \xi^1_k + \xi^2_k > 0 \} .
\end{align*}
Qualitatively, as the transmission line capacity increases,  the behavior of the two-node network begins to resemble that of a single node, driven by the aggregate net demand process $\{\xi_k^1 + \xi_k^2\}$.
That is to say, in the absence of a transmission capacity constraint between nodes 1 and 2, the nodal prices become spatially  uniform,  realizing an inter-temporal variation that depends on the nodal net demand processes only through the zero-upcrossings exhibited by their sum $\{\xi_k^1 + \xi_k^2\}$.  Intuitively, such limiting arguments reveal that an increase in a network's transmission capacity can result in an either increase, or decrease, in the LMV of storage at a node, depending on the extent to which aggregation of net demand results in nodal net demand processes exhibiting increased variability about the origin. 
}

\bibliography{Powerbib}
\appendix

\subsection{Proof of Lemma \ref{lemma:QxiLxi}}
\label{app:proof:QxiLxi}
The optimization problem in \eqref{eq:ED} is a multiparametric linear program\footnote{The objective function in \eqref{eq:ED} is piecewise linear. An optimization problem with a piecewise linear objective and linear constraints can be reformulated as a standard linear program, e.g., see \cite[Section 3.2.5]{BorrelliBook}.}, linearly parameterized in the right-hand side by $\xi$. It follows from \cite[Theorem 7.2]{BorrelliBook} that  $Q$ is continuous and piecewise affine over a polyhedral partition of $\Kcal$, which we denote by $\Lcal = \{ \Scal_1, \ldots, \Scal_{| \Lcal |} \}$. In addition, each polyhedral set $\Scal_\ell$ is full-dimensional. 
Strong duality holds in \eqref{eq:ED} for all $\xi \in \Rset^m$. Thus, the Lagrange multiplier at optimality measures the sensitivity of $Q$, if $Q$ is differentiable  \cite[Section 5.6.3]{boyd2009convex}. For each $\ell = 1,\ldots, |\Lcal |$, $Q$ is affine over $\interior{\Scal_\ell}$, and hence, differentiable. Thus, $\nabla Q(\xi) = \lambda(\xi)$ for $\xi \in \interior{\Scal_\ell}$. Since the interiors of any two distinct sets in $\Lcal$ have an empty intersection, $\Bcal = \bigcup_{\ell = 1}^{| \Lcal |} \partial\Scal_\ell$ has zero Lebesgue measure.

Suppose $v(\xi)$ is an optimizer of \eqref{eq:ED}. Then,
\begin{align}
\label{eq:subDef}
\lambda^i(\xi) \in \frac{\partial g^i (v^i)}{\partial v^i} \bigg\vert_{v(\xi)} 
=\begin{cases} 
\alpha^i, & \text{if } v^i(\xi) > 0,\\
\beta^i, & \text{if } v^i(\xi) < 0,\\
[\beta^i, \alpha^i], & \text{otherwise},
\end{cases}
\end{align}
for each $i = 1, \ldots, m$, where ${\partial g^i (v^i)} / {\partial v^i} \vert_{v(\xi)}$ denotes the  sub-differential set of $g^i$ at $v(\xi)$ with respect to the $i$-th coordinate. Then, $\alpha^i \geq \beta^i \geq 0$ implies  nonnegativity of $\lambda^i(\xi)$.


\subsection{Proof of Theorem \ref{thm:mv}}
\label{app:proof:mv}

Recall that $\Bcal = \bigcup_{\ell = 1}^{| \Lcal |} \partial\Scal_\ell$, where $\Lcal = \{ \Scal_1, \ldots, \Scal_{| \Lcal |} \}$ defines a polyhedral partition of $\Kcal$. Define 
$$\ol{\ve} := \inf_{\xi \in \Xi, s \in \Bcal} \ \  \inf_{i=1,\ldots,m} | \xi^i - s^i |.$$ 
Assumption \ref{ass:noS} guarantees that $\ol{\ve} >0$. Let $\ve \in (0,  \ol{\ve})$. Consider a single storage device of capacity $\ve$ installed at node $i$, i.e., $b = \ve \bone^i$. Here,  $\bone^i$ denotes the $i^{\rm th}$ standard basis vector of appropriate dimension. 
To prove \eqref{eq:mv}, it suffices to establish the following identity 
\begin{align}
J^*(0)  -  J^*(\ve \bone^i) &= \ve \E \left[ \sum_{k=0}^{N-2}({\lambda}^i_{k+1 | k} - {\lambda}^i_k)^+ \right].\label{eq:mvGoal.1}
\end{align}
We establish the desired form of the optimal cost through an argument based on dynamic programming (DP). For each $z \in \Zcal(b)$ and $\xi_{\leq k} \in \Xi^{k+1}$, define the optimal value functions:
\begin{align}
\begin{aligned}\label{eq:JN1}
J_{N-1}^*(z, \xi_{\leq N-1}; b)   :=  & \underset{u \in \Rset^m, v \in \Rset^m}{\text{minimum}}  && g(v), \\
&  \ \ \text{subject to} &&  z + u \in \Zcal (b), \\
& \quad && v - \xi_{N-1} - u \in \Pcal
\end{aligned}
\end{align}
and
\begin{align} 
\begin{aligned}
& J_{k}^*(z, \xi_{\leq k}; b) \\
& \quad  := \underset{u \in \Rset^m, v \in \Rset^m}{\text{minimum}}  && g(v) + \E\left[ J_{k+1}^*(z+u, \xi_{\leq k+1}; b) \ | \ \xi_{\leq k} \right], \\
& \quad \quad \ \   \text{subject to} && z + u \in \Zcal (b), \\
& \quad  &&  v - \xi_{k} - u \in \Pcal
\end{aligned}\label{eq:Jk}
\end{align}
for $k=0, \ldots, N-2$. By \cite[Proposition 1.3.1]{bertsekas1995dynamic}, the policy obtained as a recursive solution to the above system of DP equations is indeed optimal, and 
$$ J^*(b) = \E \left[ J_0^*(z_0 = 0, \xi_0; b) \right]$$ 
for any $b \in \Rset^m_+$. To emphasize the optimal policy's dependence on the storage capacity parameter $b$, we write $\mu_k^*(z_{\leq k}, \xi_{\leq k}; b)$ and $\nu_k^*(z_{\leq k}, \xi_{\leq k}; b)$  for $k=0, \ldots, N-1$. 
We have the following result, from which equation \eqref{eq:mvGoal.1}, and hence our desired result  \eqref{eq:mv} -- is an immediate consequence. We defer its proof to Appendix \ref{app:lemma:JkVk}.

\begin{lemma}
\label{lemma:JkVk}
For each $\ve \in (0, \ol{\ve})$, $z \in \Zcal(\ve \bone^i)$, and $\xi_{\leq N-1} \in \Xi^N$,
\begin{align}
& J_{k}^*(0, \xi_{\leq k}; 0) - J_{k}^*(z, \xi_{\leq k}; \ve \bone^i) \notag\\
& \qquad =  \lambda^i_k z^i + \ve \E \left[ \sum_{j=k}^{N-2} \left( {\lambda}^i_{j+1 | j}- \lambda^i_j \right)^+ \ \bigg\vert \ \xi_{\leq k}\right], \label{eq:JkVk}
\\
& \mu^*_{k}(z_{\leq k}, \xi_{\leq k}; \ve \bone^i ) = \begin{cases} \left( \ve - z^i_{k} \right) \bone^i, & \text{if } \lambda^i_{k} \leq {\lambda}^i_{k+1 | k},\\ 
- z^i_{k} \bone^i, & \text{otherwise}
\end{cases}
\label{eq:muk}
\end{align}
for $k=0, \ldots, N-2$, and 
\begin{align}
& J_{N-1}^*(0, \xi_{\leq N-1}; 0) - J_{N-1}^*(z, \xi_{\leq N-1}; \ve \bone^i) =  \lambda^i_{N-1} z^i, \label{eq:JN1VN1}
\\
& \mu^*_{N-1}(z_{\leq N-1}, \xi_{\leq N-1}; \ve \bone^i) = -z^i_{N-1} \bone^i.
\label{eq:muN1}
\end{align}
\end{lemma}
To establish the inequality \eqref{eq:mvUp}, we have the following string of arguments.
\begin{align}
\E \left[ \left( {\lambda}^i_{k+1 | k}- \lambda^i_k \right)^+  \right]
&=   \E \left[ \left(  \E \left[{\lambda}^i_{k+1} - \lambda^i_k \ \vert \ \xi_{\leq k} \right]  \right)^+ \right] \notag\\
&\leq  \E \left[ \E \left[  \left( {\lambda}^i_{k+1} - \lambda^i_k   \right)^+ \vert \ \xi_{\leq k} \right] \right] \label{eq:ineqG.1}\\
&=  \E \left[  \left( {\lambda}^i_{k+1} - \lambda^i_k   \right)^+ \right], \label{eq:ineqG.2}
\end{align}
where \eqref{eq:ineqG.1} follows from Jensen's inequality on the convex function $F(x) = (x)^+$ for $x \in \Rset$, and \eqref{eq:ineqG.2} follows from the law of iterated expectation.
Now, $(x)^+ = \frac{1}{2}( | x | + x )$ for any $x \in \Rset$, and hence,

\begin{align*}
\E \left[ \left( {\lambda}^i_{k+1}- \lambda^i_k \right)^+  \right]
&=  \frac{1}{2}\E \left[ \left\vert {\lambda}^i_{k+1}  - \lambda^i_k\right\vert \right] + \frac{1}{2} \E\left[ {\lambda}^i_{k+1}  - \lambda^i_k \right].
\end{align*}
Plugging the above expression into \eqref{eq:ineqG.2}, and summing both sides from $k=0$ to $k=N-2$ yields the upper bound in \eqref{eq:mvUp}.

\subsection{Proof of Lemma \ref{lemma:JkVk}}
\label{app:lemma:JkVk}
The proof proceeds by backward induction. For $z \in \Zcal(\ve\bone^i)$, \eqref{eq:JN1} can be written as 
\begin{align}
& J_{N-1}^*(z, \xi_{\leq N-1}; \ve \bone^i)  \notag \\
& \qquad = \underset{\substack{u^i \in \Rset, \\ 0 \leq z^i + u^i  \leq \ve}}{\text{minimum}} \ \   \underset{\substack{v \in \Rset^m, \\
v - \xi_{N-1} - u^i \bone^i \in \Pcal}}{\text{minimum}} \quad g(v) \notag \\
&\qquad  =  \underset{\substack{u^i \in \Rset, \\ 0 \leq z^i + u^i  \leq \ve}}{\text{minimum}} \ \ Q( \xi_{N-1} + u^i \bone^i),\label{eq:JN1.Qrep}
\end{align}
where $Q$ is defined as in \eqref{eq:ED}. 
Assumption \ref{ass:noS} guarantees that $\xi_{N-1} \in \interior{\Scal}$, for some $\Scal \in \Lcal$. Also, our choice of $\ve  \in (0, \ol{\ve})$ implies that $\xi_{N-1} + u^i \bone^i \in \interior{\Scal}$ for the same $\Scal \in \Lcal$. Now, by  Lemma \ref{lemma:QxiLxi}, we have that $Q$ is affine and  $\nabla Q(\xi) = \lambda(\xi)$ for all $\xi \in \interior{\Scal}$. Hence,
\begin{align}
Q( \xi_{N-1} + u^i \bone^i) 
& = Q(\xi_{N-1}) + {\lambda}^i_{N-1} u^i.
\label{eq:Q}
\end{align}
Plugging the above expression into \eqref{eq:JN1.Qrep}, we get
\begin{align}
J_{N-1}^*(z, \xi_{\leq N-1}; \ve \bone^i) = {Q(\xi_{N-1})} + \underset{\substack{u^i \in \Rset, \\ 0 \leq z^i + u^i  \leq \ve}}{\text{minimum}} \ \   {\lambda}^i_{N-1} u^i.
\label{eq:JN1.1}
\end{align}

Now, $Q(\xi_{N-1}) = J_{N-1}^*(0, \xi_{\leq N-1}; 0)$. Also, $\lambda^i_{N-1} \geq 0$ implies that $(u^i)^* = -z^i$ is the optimizer in \eqref{eq:JN1.1}. That in turn proves \eqref{eq:muN1}, and we obtain 
\begin{align*}
J_{N-1}^*(z, \xi_{\leq N-1}; \ve \bone^i) 
&= J_{N-1}^*(0, \xi_{\leq N-1}; 0) - {\lambda}^i_{N-1} z^i.
\end{align*}
The above equation is precisely \eqref{eq:JN1VN1}.

Continuing the induction hypothesis, suppose that \eqref{eq:JkVk} -- \eqref{eq:muk} are satisfied for periods $k+1, \ldots, N-2$. In what follows, we prove \eqref{eq:JkVk} -- \eqref{eq:muk} for period $k$. To that end, we first rewrite \eqref{eq:Jk} as
\begin{align*}
J_{k}^*(z, \xi_{\leq k}; \ve \bone^i) 
& = \underset{\substack{u^i \in \Rset, \\ 0 \leq z^i + u^i  \leq \ve}}{\text{minimum}} 
\ \ f(u^i)
\end{align*}
for each $z \in \Zcal(\ve \bone^i)$. 
Here,  $f(u^i)$ is defined as
\begin{align}
f(u^i) 
& := \E \left[ J^*_{k+1}(z + u^i \bone^i, \xi_{\leq k+1}; \ve \bone^i ) \ | \ \xi_{\leq k}\right] \notag \\
& \qquad \quad +  \underset{\substack{v \in \Rset^m, \\
v - \xi_{k} - u^i \bone^i \in \Pcal}}{\text{minimum}} \ \ g(v).
\label{eq:defFu}
\end{align}
The above expression for $f(u^i)$ contains two terms. We tackle them individually. Using the induction hypothesis and the law of iterated expectation, the first term in \eqref{eq:defFu} satisfies
\begin{align*}
& \E\left[ J_{k+1}^*(z + u^i \bone^i, \xi_{\leq k+1}; \ve \bone^i) \ | \ \xi_{\leq k}\right]  \\
& \ =  \E\left[ J_{k+1}^*(0, \xi_{\leq k+1}; 0) \ | \ \xi_{\leq k} \right] - \E \left[ \lambda^i_{k+1} (z^i + u^i) \ \vert \ \xi_{\leq k} \right] \\
& \ \ \ \ - \ve \E \left[ \E \left[ \sum_{j=k+1}^{N-2} \left( {\lambda}^i_{j+1 | j} - \lambda^i_j \right)^+   \bigg\vert \ \xi_{\leq k+1} \right]  \bigg\vert \ \xi_{\leq k} \right] \\
& \ = \E\left[ J_{k+1}^*(0, \xi_{\leq k+1}; 0) \ | \ \xi_{\leq k} \right] -  \lambda^i_{k+1 | k} u^i \\
& \ \ \ \ -\lambda^i_{k+1 | k} z^i - \ve \E \left[ \sum_{j=k+1}^{N-2} \left( {\lambda}^i_{j+1 | j} - \lambda^i_j \right)^+   \bigg\vert \ \xi_{\leq k} \right].
\end{align*}
And, the second term in \eqref{eq:defFu} is precisely $Q(\xi_{k} + u^i \bone^i)$, which further satisfies
$$Q(\xi_{k} + u^i \bone^i) = Q(\xi_k) + {\lambda}^i_k  u^i.$$
Plugging the two derived expressions into \eqref{eq:defFu}, we obtain
\begin{align}
&J_{k}^*(z, \xi_{\leq k}; \ve \bone^i) \notag\\
& \  = Q(\xi_k) + \E\left[ J_{k+1}^*(0, \xi_{\leq k+1}; 0) \ | \ \xi_{\leq k} \right]  \notag \\
& \ \quad +  \underset{\substack{u^i \in \Rset, \\ 0 \leq z^i + u^i  \leq \ve}}{\text{minimum}} 
\left( {\lambda}^i_k  -  \lambda^i_{k+1 | k} \right) u^i \notag\\
&  \ \quad -  \lambda^i_{k+1 | k} z^i - \ve \E \left[ \sum_{j=k+1}^{N-2} \left( {\lambda}^i_{j+1 | j} - \lambda^i_j \right)^+   \bigg\vert \ \xi_{\leq k} \right].
\label{eq:Jk.5}
\end{align}
In order to prove \eqref{eq:JkVk} -- \eqref{eq:muk}, we need to further simplify \eqref{eq:Jk.5}. The first two terms in the right-hand side of \eqref{eq:Jk.5} simplify to
\begin{align}
Q(\xi_k) + \E \left[ J_{k+1}^*(0, \xi_{\leq k+1}; 0) \ | \ \xi_{ \leq k}  \right] = J_{k}^*(0, \xi_{\leq k}; 0).
\label{eq:Jk0}
\end{align}
And, the third term in \eqref{eq:Jk.5} can be written as
\begin{align}
& \underset{\substack{u^i \in \Rset, \\ 0 \leq z^i + u^i  \leq \ve}}{\text{minimum}} 
\left( {\lambda}^i_k  -  \lambda^i_{k+1 | k} \right) u^i  \notag \\
& \qquad = \left(  \lambda^i_{k+1 | k} - \lambda^i_k \right) z^i + \ve \left(  \lambda^i_{k+1 | k} - \lambda^i_k \right)^+. \label{eq:Jk.opt}
\end{align}
This follows from noting that $(u^i)^* = \ve - z^i$ is the optimizer, if $\lambda^i_k \leq \lambda^i_{k+1 | k}$; otherwise, the optimizer is $(u^i)^* = -z^i$. This proves \eqref{eq:muk}. Finally, \eqref{eq:JkVk} follows from combining the derived expressions in \eqref{eq:Jk0} and \eqref{eq:Jk.opt} back into \eqref{eq:Jk.5}, and simplifying. The details are omitted for brevity.


\subsection{Proof of Proposition \ref{prop:spCase}}
\label{app:proof:spCase}
In this proof, we only establish the result that $\lambda^i(\xi) \in \{ \alpha, \beta \}$ for all $\xi \in \Xi$. The remaining results are an immediate consequence. We do so by showing that $\Fcal \subseteq \Bcal$, where
\begin{align}
\label{eq:spCaseGoal} \Fcal := \{{\xi} \in \Kcal \ \vert \ \lambda^i (\xi) \notin \{ \alpha, \beta \} \}.
\end{align}
The following alternative description of the feasible power injection polytope $\Pcal$ defined in \eqref{eq:injPoly} shall prove useful in the sequel.
$$\Pcal = \left\{ x \in \Rset^m \ \vert \ -f  \leq H x \leq  f, \ \bone^\top x = 0 \right\},$$ 
where $H = B \left( Y^\top Y + \bone^1 (\bone^1)^\top \right)^{-1} Y^\top$, and $\bone$ denotes a vector of all ones of appropriate size.

Let $\xi \in \Fcal$ and $v(\xi)$ be an optimizer of \eqref{eq:ED} with net demand $\xi$. Then, \eqref{eq:subDef} implies that $\lambda^i(\xi)$ belongs to the set of sub-differentials of $g^i$ with respect to the $i$-th coordinate at $v(\xi)$, and hence,
\begin{align}
\label{eq:vi0}
\xi \in \Fcal \implies \lambda^i(\xi) \notin  \{ \alpha, \beta \} \implies v^i (\xi) = 0.
\end{align}
Then, $H(v(\xi) - \xi) \in \Rset^\ell$ defines the vector of power flows on the $\ell$ transmission lines at the optimum of \eqref{eq:ED}. One of two cases can arise: (1) the power flow on each transmission line has a magnitude strictly less than its capacity, or (2) the power flow on at least one transmission line equals its capacity. Each case is analyzed separately.

\textit{Case 1:} We argue that $\bone^\top \xi = 0$ in this case. To that end, suppose $\bone^\top \xi \neq 0$ to the contrary. Since, $\bone^\top v(\xi) = \bone^\top \xi \neq 0$, there exists a node $j \neq i$ in the network, for which $v^j(\xi) \cdot (\bone^\top \xi) > 0$. For $\delta > 0$, define
$$v_\delta (\xi) := v(\xi) + \delta (\bone^\top \xi) (\bone^i - \bone^j).$$
By hypothesis, $ f- | H(v(\xi)-\xi)| > 0$, where the inequality is element-wise, and $| \cdot |$ is the absolute value operator. One can then choose $\delta> 0$ small enough to satisfy
$$\delta | H(\bone^i - \bone^j) | \leq f- | H(v(\xi) - \xi) | \quad \text{and} \quad \delta < | v^j (\xi) |.$$
For such a $\delta$, we have $v_\delta (\xi) - \xi \in \Pcal$, and $g(v(\xi)) = g(v_\delta (\xi))$. Hence, $v_\delta (\xi)$ is another optimizer of \eqref{eq:ED}. However, $v^i_\delta (\xi) \neq 0$ violates \eqref{eq:vi0}, resulting in a contradiction.

\textit{Case 2:} At the optimum, let the power flows on certain transmission lines equal the respective line capacities. We cut the graph of the power network at these transmission lines. The network being acyclic, said cut results in a \emph{forest} of disconnected acyclic graphs. Consider the connected component containing bus $i$; call it $\Ccal(i)$. Let $\xi^{\Ccal(i)}$, 
$v^{\Ccal(i)}(\xi)$, $Y^{\Ccal(i)}$, $B^{\Ccal(i)}$ and $f^{\Ccal(i)}$ be the corresponding vectors and matrices with rows and/or columns restricted to the nodes in $\Ccal(i)$. Define $\Ccal'(i) := \{ 1, \ldots, m \} \setminus \Ccal(i)$. At the optimum, the power flow on any transmission line joining two neighboring nodes in $\Ccal(i)$ is not its capacity, but between a node in $\Ccal(i)$ and its neighboring node in $\Ccal'(i)$ is at its capacity. Let $a \in \Ccal(i)$ and $a' \in \Ncal(a) \subseteq \Ccal'(i)$, where $\Ncal(a)$ is the set of neighbors of $a$ in $\Ccal'(i)$. The power flow $p^{a a'}$ from $a$ to $a'$ satisfies $p^{a a'} = \pm f^{a a'}$, where $f^{a a'}$ is the capacity of the corresponding line. Then, $v^{\Ccal(i)}(\xi)$ solves an optimization problem similar to \eqref{eq:ED}, restricted to the nodes in $\Ccal(i)$ with the net demand vector
\begin{align}
\label{eq:xiTilde}
\tilde{\xi}^{\Ccal(i)} := \xi^{\Ccal(i)} + \sum_{a \in \Ccal(i), a' \in \Ncal(a) } p^{a a' } \bone^a.
\end{align}
Arguing similar to case 1, one can show that $\bone^\top \tilde{\xi}^{\Ccal(i)} = 0$.

To finally show $\Fcal \subseteq \Bcal$, define $\tilde{\xi}^{\Ccal}$ as in \eqref{eq:xiTilde} for any connected component $\Ccal$ of the power network. Consider the union of the sets
$ \{ \xi \in \Kcal \ | \ \bone^\top \tilde{\xi}^{\Ccal}=0 \}$
over all connected components $\Ccal$. Said set is a finite union of hyperplanes in $\Rset^m$ of dimension less than $m$, and hence, has zero Lebesgue measure. Moreover, this set contains $\Fcal$. Hence, $\Fcal$ has zero Lebesgue measure. Suppose $\Fcal \cap \interior{\Scal}$ is nonempty for some set $\Scal \in \Lcal$. Lemma \ref{lemma:QxiLxi} then implies $\lambda(\xi)$ is constant over $\interior{\Scal}$, and hence, $\interior{\Scal} \subseteq \Fcal$. However, $\interior{\Scal}$ is full-dimensional and has positive Lebesgue measure, while $\Fcal$ has zero measure. This yields a  contradiction, implying $\Fcal \subseteq \Bcal$.


\end{document}